\documentclass[reqno,12pt]{amsart}

\usepackage{amsmath,amsthm,amsfonts,amssymb}

\date{}
\begin{document}
{\large

\centerline{{\sc On topological spaces possessing
}}

\vskip .1in

\centerline{{\sc
 uniformly distributed sequences}}

\vskip .2in

\centerline{{\sc Bogachev V.I., Lukintsova M.N.}}

\vskip .3in
\centerline{{\bf Abstract}}

\vskip .1in

{\small Two classes of topological spaces are introduced
on which every probability Radon  measure
possesses a uniformly distributed  sequence
or a uniformly tight uniformly distributed sequence.
It is shown that these classes are stable under multiplication by
completely regular Souslin spaces.}

\vskip .2in

We introduce two classes of topological spaces,
on which every probability Radon  measure
possesses a uniformly distributed  sequence
or a uniformly tight uniformly distributed sequence.
It is shown that these properties are preserved under multiplication by
completely regular Souslin spaces. The concepts related to measures
on topological spaces and weak
convergence of measures which we use below can be found in~[1].

Let $X$ be a completely regular topological space and let
$C_b(X)$ be the space of bounded continuous functions on~$X$.
We recall that a nonnegative measure $\mu$ on the Borel
$\sigma$-algebra of the space $X$ is called Radon
if, for every Borel set $B$ and every
$\varepsilon>0$, there exists a  compact set $K\subset B$ such that
$\mu(B\backslash K)<\varepsilon$. A~sequence of
nonnegative Radon measures $\mu_n$ on $X$ is called
uniformly tight if, for every $\varepsilon>0$,
there exists a compact set $K_\varepsilon\subset X$ such that
$\mu_n(X\backslash K_\varepsilon)<\varepsilon$  for all~$n$.
A~sequence of Radon measures $\mu_n$  converges weakly to a
Radon measure $\mu$ if the integrals of every bounded continuous function
with respect to the measures $\mu_n$ converge to the integral of this
function with respect to the measure~$\mu$.
A~sequence of points $x_n\in X$ is called
uniformly distributed  with respect to a probability Radon measure  $\mu$
on~$X$ if, for every bounded continuous
function $f$ on~$X$, one  has the equality
$$
\int_X f(x)\, \mu(dx)=\lim\limits_{n\to\infty}
\frac{f(x_1)+\cdots+f(x_n)}{n}.
$$
Let $\delta_a$ denote the probability measure concentrated
at the point $a\in X$. The definition means weak convergence
 of the  measures $n^{-1}(\delta_{x_1}+\cdots+\delta_{x_n})$ to the measure~$\mu$.
According to Niederreiter's theorem [2] (see also
[1], \S8.10(ix)), the existence of a uniformly
distributed sequence with respect to $\mu$
is equivalent to the fact that some sequence of probability measures with finite
supports converges weakly to~$\mu$.
If $X$ is a completely regular Souslin space,
then every Radon probability measure on $X$ has a
uniformly distributed sequence. In general, this is not the case even for compact
spaces: the Stone--\v{C}ech  com\-pac\-ti\-fi\-ca\-tion of the space of natural
numbers serves as a counterexample
 (see [1], Example 8.10.54).
Certainly, the notion of a uniformly distributed sequence is meaningful also for
Baire measures.

\vskip .1in {\bf Definition 1.}
{\it We shall say that a Radon probability
 measure $\mu$ on $X$ has  a $t$-uniformly distributed sequence if there exists
 a sequence of points $x_n\in X$ such that the
sequence of measures $n^{-1}(\delta_{x_1}+\cdots+\delta_{x_n})$
is uniformly tight and converges  weakly to~$\mu$.}

\vskip .1in
{\bf Definition 2.}
{\it
We shall say that $X$ possesses  property {\rm(ud)}
if every Radon probability measure  on $X$ has a uniformly distributed sequence.
If every such a measure has a
$t$-uniformly distributed sequence, then we call $X$ a  space with property {\rm(tud)}.
}
\vskip .1in

By the above cited result of Niederreiter property  (ud) is equivalent to
the property that every  Radon measure on $X$ is a weak limit of a sequence
of probability measures with finite supports, and,
as one can see from the proof of Niederreiter's theorem (see [1], Theorem 8.10.52),
 property (tud) requires additionally the uniform tightness of some
such sequence.

We note that the indicated properties are not always preserved by passing to a
 compact subset. For example, under the continuum hypothesis,
the product of the  continuum of intervals has property~(ud), but this product contains
a closed subset homeomorphic to the Stone--\v{C}ech compactification of the space of
natural numbers, which  has no property~(ud). In particular, supports
of approximating discrete  measures cannot be always chosen in the support
of the measure which is approximated.
Diverse results related to uniformly distributed sequences in topological
space can be found in [3]--[7].

The following lemma shows, in particular,
 that it suffices to verify the in\-tro\-du\-ced properties only for
measures with compact supports.

\vskip .1in
{\bf Lemma.}
{\it
Suppose that for every $n\in\mathbb{N}$ a subspace $X_n$ in $X$ is
 measurable with respect to all Radon measures on $X$ and has property~{\rm(tud)}.
Then the space
$Y:=\bigcup_{n=1}^\infty X_n$ possesses this property as well.
The analogous assertion is true for property~{\rm(ud)}.
}
\vskip .1in
{\bf Proof.}
Let $\mu$ be a Radon probability measure  on~$Y$. Set
$X_0:=\emptyset$.
Let $I_A$ denote the indicator function of the set~$A$
and let  $I_A\cdot\mu$ denote the measure with density $I_A$ with respect to~$\mu$.
Properties (tud) and (ud) are preserved under finite unions.
Indeed, if $\mu$ is a Radon probability  measure on $X_1\cup X_2$, then one can find
nonnegative discrete measures $\mu_n^1$ on $X_1$ with
$\mu_n^1(X_1)=\mu(X_1)$ which are weakly convergent on $X_1$
to the measure $I_{X_1}\cdot\mu$, and also
nonnegative discrete measures $\mu_n^2$ on $X_2$ with
$\mu_n^2(X_2)=\mu(X_2\backslash X_1)$ which are weakly convergent on $X_2$
to the measure $I_{X_2\backslash X_1}\cdot\mu$. Then the discrete probability
measures $\mu_n^1+\mu_n^2$  converge weakly to $\mu$ on $X_1\cup X_2$.

Hence we may assume further that $X_n\subset X_{n+1}$.
By hypothesis, for every $n$ there exists a uniformly tight sequence of
nonnegative measures $\mu_{n,k}$ on $X_n$ with finite supports and
$\|\mu_{n,k}\|=\mu(X_n\backslash X_{n-1})$ which
converges weakly to the measure $I_{X_n\backslash X_{n-1}}\cdot\mu$ on $X_n$
as $k\to\infty$. Set
$$
\nu_n:=c_n^{-1}[\mu_{1,n}+\cdots+\mu_{n,n}],
$$
$$
c_n:=\mu_{1,n}(X)+\cdots+\mu_{n,n}(X)=\mu(X_n).
$$
 One obtains probability measures with finite supports in~$Y$. It is clear that $c_n\to 1$
as $n\to\infty$. The measures $\nu_n$  converge weakly to $\mu$ on~$Y$.
Indeed, let $f$ be a bounded continuous function on $Y$ and let $\varepsilon>0$.
We may assume that $|f|\le 1$. Let us find $n_1$ such that
$$
\sum_{n=n_1}^\infty \mu(X_n\backslash X_{n-1})<\varepsilon
$$
and
$$
\|\nu_n- \mu_{1,n}-\cdots-\mu_{n,n}\|\le \varepsilon,
\quad \forall\, n\ge n_1.
$$
Next we find $m\ge n_1$ such that
$$
\Bigl|\int_{X_n\backslash X_{n-1}} f\, d\mu
- \int_{X_n} f\, d\mu_{n,k}\Bigr|<\varepsilon n_1^{-1}
$$
for every $n=1,\ldots,n_1$ and all $k\ge m$.
Then, for all $n\ge m$ we obtain
$$
\Bigl|\int_{Y} f\, d\mu - \int_{Y} f\, d\nu_n\Bigr|
\le 4\varepsilon .
$$
It remains to observe that the sequence $\{\nu_n\}$ is
uniformly tight on~$Y$. To this end, for fixed $\varepsilon>0$ we
choose $n_1$ as above and then for every
$n=1,\ldots,n_1$ we find a compact set $K_n\subset X_n$ such that
$\mu_{n,k}(X\backslash K_n)\le \varepsilon$ for all~$k$.
The set $K=K_1\cup\cdots\cup K_{n_1}$ is compact and
$\nu_n(Y\backslash K)\le 3\varepsilon$ for all~$n$.
The case of property (ud) is similar.

\vskip .1in
{\bf Remark 1.}
It follows from this lemma that every  probability Radon
con\-cen\-tra\-ted
on a countable union of metrizable compact sets has a $t$-uniformly
distributed sequence that belongs to this union. Hence any
completely regular
space in which all compact sets are metrizable has property (tud).
For example, this is the case for completely regular Souslin spaces.

\vskip .1in
{\bf Theorem.}
{\it Suppose that a space $X$ possesses property {\rm(tud)}.
Then, the  space
$X\times Y$ has this property for every nonempty completely regular
space $Y$ in which all compact sets are metrizable.
The analogous assertion is true for property~{\rm(ud)}.}
\vskip .1in

{\bf Proof.}
Since the projections of every measure on $X\times Y$ onto the factors are concentrated
on countable unions of compact sets, and all compact sets in $Y$ are metrizable,
the lemma reduces the general case to the case where $X$ and $Y$ are compact
and $Y$ is a  metric space. In this case properties  (ud)  and (tud) coincide.

Let $\mu$ be a Radon probability measure on $X\times Y$, let $\nu$ be its
projection onto~$X$, and let $\mu_Y$ be its projection onto~$Y$.
It is known see [1], \S10.4) that there exist
conditional Radon probability  measures $\mu^x$ on~$Y$, i.e.,
$$
\int_{X\times Y} f(x,y)\, \mu(dxdy)
=\int_X\int_Y f(x,y)\, \mu^{x}(dy)\, \nu(dx)
$$
for every bounded Borel function $f$ on $X\times Y$,
and for every Borel set $B\subset Y$ the function
$x\mapsto \mu^x(B)$ is measurable with respect to~$\nu$.
Let us equip the space $\mathcal{P}(Y)$ of  Borel probability measures on $Y$ by the
Kantorovich--Rubinshtein metric $d$ (see [1], \S8.3),
which is defined by the formula
$$
d(\mu_1,\mu_2)=\sup \Bigl|
\int_Y f(y)\, \mu_1(dy)- \int_Y f(y)\, \mu_2(dy)\Bigr|,
$$
where sup is taken over all functions $f$ on $Y$ such that
$|f(y)|\le 1$ and $f$ is Lipschitzian with constant~$1$.
 Then $(\mathcal{P}(Y),d)$ is a compact metric space (see [1],
 Theorem 8.3.2 and Theorem~8.9.3). It is easily verified that the mapping
$\xi\colon\, x\mapsto \mu^x$ from $X$ to $\mathcal{P}(Y)$ is measurable with
respect to the measure~$\nu$. To this end we observe that there is a countable set
in $C(Y)$ separating the measures on~$Y$, and for every function
$\varphi\in C(Y)$ the function
$$
x\mapsto \int_Y \varphi(y)\, \mu^x(dy)
$$
is measurable with respect to~$\nu$. By Lusin's theorem there exists a sequence
of  compact sets $K_n\subset X$ such that
$\nu(K_n)>1-2^{-n}$, $K_n\subset K_{n+1}$ and
the restrictions of the mapping $\xi$ to $K_n$ are continuous.
We recall that any continuous mapping $G$ on a compact $K$ in the space $X$
taking values in a Banach space $E$ extends to a continuous mapping $G_K$ from $X$
to the closed convex hull of the set~$G(K)$. The completeness of $E$ can be omitted
if the closed convex hull of $G(K)$ is complete (it suffices to consider the completion
of~$E$). Since $\mathcal{P}(Y)$ is a complete convex subset of the normed
space of all Radon measures on $Y$ with the Kantorovich--Rubinshtein norm
(this normed space is not Banach if the compact $Y$ is infinite),
we obtain that for every $n$ there exists a continuous mapping
$\xi_n\colon\, x\mapsto \xi_n^x$, $X\to \mathcal{P}(Y)$, that coincides with
 $\xi$ on~$K_n$. Set $K_0=\emptyset$.

For every fixed $n$, for each $i=1,\ldots,n-1$ we find an open set $U_{n,i}\supset K_i$
with the following properties:
$$
d(\xi_i^x,\xi_n^x)\le 2^{-n}
\quad\hbox{for all
$i=1,\ldots,n-1$ and all $x\in U_{n,i}$.}
$$
 This is possible due to the coincidence of $\xi_n$ and $\xi_i$ on the compact set
 $K_i$ and the continuity of  $\xi_n$ and~$\xi_i$. We may assume
 that $U_{n,i}\subset U_{n,i+1}$ whenever $i\le n-2$
by dealing with the sets $U_{n,1}\cap\cdots\cap U_{n,i}$, since $K_i\subset K_{i+1}$.

By our assumption, for every fixed $m$ there is a sequence of nonnegative
Radon measures $\nu_{m,i}$ with finite supports $S_{m,i}$ and
$\nu_{m,i}(X)=\nu(K_m\backslash K_{m-1})$ that
converge weakly on $X$ to the measure~$I_{K_m\backslash K_{m-1}}\cdot\nu$.
Now for every $n$ we find a number $m_n\ge n$ such that
$$
\nu_{1,m_n}(X\backslash U_{n,1})
+
\nu_{2,m_n}(X\backslash U_{n,2})
+\cdots
+
\nu_{n-1,m_n}(X\backslash U_{n,n-1})
\le 2^{-n}. \eqno(1)
$$
This is possible because weak convergence of the measures
$\nu_{l,i}$ to the measure $I_{K_l\backslash K_{l-1}}\cdot\nu$
gives the relationship
$\limsup\limits_{i\to\infty}\nu_{l,i}(Z)\le
\nu(Z\cap (K_l\backslash K_{l-1}))$ for all closed
sets~$Z$, and we have
$$
\nu((X\backslash U_{n,l})\cap (K_l\backslash K_{l-1}))=0
$$
by the inclusion $K_l\subset U_{n,l}$ for all $l\le n-1$.
Set
$$
\nu_n:=\nu_{1,m_n}+\cdots +\nu_{n-1,m_n}.
$$
Let us observe that $\nu_n$ is a nonnegative measure with finite
support and
$$
\nu_n(X)=\nu_{1,m_n}(X)+\cdots +\nu_{n-1,m_n}(X)=\nu(K_{n-1}).
$$
Therefore, the  measures $\nu_n(X)^{-1}\nu_n$  are probability measures
for all $n>1$.

Let us show that the measures $\nu_n$  converge weakly to~$\nu$.
Let $\psi\in C_b(X)$, $|\psi(x)|\le 1$ and
$\varepsilon>0$. Let us take $m$ with $2^{-m}<\varepsilon$.
By  weak convergence of the measures $\nu_{l,i}$ to the measures
$I_{K_l\backslash K_{l-1}}\cdot\nu$ for each fixed  $l=1,\ldots,m$,
there exists $N>m$ such that for all $l=1,\ldots,m$ and $i\ge N$
one has the inequality
$$
\biggl|
\int_X \psi(x)\, \nu_{l,i}(dx)
-\int_{K_l\backslash K_{l-1}} \psi(x)\, \nu(dx)
\biggr|\le \varepsilon m^{-1}.
$$
Whenever $i\ge N$ we obtain
$$
\biggl|
\int_X \psi(x)\, (\nu_{1,i}+\cdots+\nu_{m,i})(dx)
-\int_{K_m} \psi(x)\, \nu(dx)
\biggr|\le \varepsilon .
$$
For all $n\ge N$ this yields the estimate
$$
\biggl|
\int_X \psi(x)\, \nu_{n}(dx)
-\int_{X} \psi(x)\, \nu(dx)
\biggr|\le 3\varepsilon
$$
by the relationships
$\nu(X\backslash K_m)\le \varepsilon$
and
$$
\|\nu_n-\nu_{1,m_n}-\cdots-\nu_{m,m_n}\|
\le
\nu(K_{m+1}\backslash K_m)
+\nu(K_{m+2}\backslash K_{m+1})+\cdots \le
2^{-m}\le \varepsilon.
$$
Thus, weak convergence of $\nu_n$ to $\nu$ is established.

By the compactness of $(\mathcal{P}(Y),d)$, for every fixed $n$ one can divide
$X$ into pairwise disjoint Borel parts $B_{n,i}$, $i=1,\ldots,N_n$, such
that $d(\xi_n^x,\xi_n^z)<2^{-n}$ if $x,z\in B_{n,i}$.
Let us choose a point $x_{n,i}$ in $B_{n,i}$ and find a
measure $\sigma_{n,i}\in\mathcal{P}(Y)$ with finite support such
that $d(\xi_n^{x},\sigma_{n,i})\le 2^{1-n}$ whenever $x\in B_{n,i}$.
For this purpose it suffices to  take such a measure in the ball of radius $2^{-n}$
centered at~$\xi_n^{x_{n,i}}$.
Set
$$
\mu_n^x:=\sigma_{n,i}\quad\hbox{if $x\in B_{n,i}$.}
$$
Therefore,
$$
\sup_{x\in X} d(\xi_n^x,\mu_n^x)\le 2^{1-n}.
$$
Finally, we define a measure $\mu_n$ on $X\times Y$ with finite support as
 the measure with the projection $\nu_n$ on $X$ and the conditional
measures $\mu_n^x$, i.e.,
$$
\mu_n(B)=\int_X\int_Y I_B(x,y)\, \mu_n^x(dy)\, \nu_n(dx).
$$
Let us show that the measures $\mu_n$  converge weakly
to~$\mu$. It suffices to establish convergence on the functions of the form
$$
f(x,y)=\psi(x)\varphi(y),
\quad
\hbox{where $\psi\in C_b(X)$,
$\sup_x |\psi(x)|\le 1$,
$\sup_y |\varphi(y)|\le 1$,}
$$
and the function $\varphi$ is Lipschitzian with constant~$1$.
Indeed, it is readily seen that the linear span of the set of such functions
is dense in $C(X\times Y)$.
 Let $\varepsilon >0$. We choose $m$ such that
 $2^{-m}<\varepsilon$.
There holds the inequality
$$
\biggl|
\int_{X}\psi(x)\int_Y
\varphi(y)\, \xi_m^x (dy) \, \nu(dx)-
\int_{X}
\psi(x)
\int_Y
\varphi(y)\, \xi^x (dy) \, \nu(dx)
\biggr|
\le 2^{-m} \eqno(2)
$$
since $\xi_m=\xi$ on $K_m$ and
$\nu(X\backslash K_m)\le 2^{-m}$.
We show that for all $n\ge m$ one has the inequality
$$
\biggl|
\int_X \psi(x)
\int_Y \varphi(y)\, \xi_m^x (dy)\, \nu_n(dx)
-
\int_X \psi(x)
\int_Y \varphi(y)\, \xi_n^x (dy)\, \nu_n(dx)
\biggr|\le 2^{2-m}. \eqno(3)
$$
Indeed, the  measure $\nu_n$ differs in the variation norm from the measure
$$
\eta:=\nu_{1,m_n}+\cdots +\nu_{m,m_n}
$$
not greater than in~$2^{-m}$.
By estimate (1) and the inclusion $U_{n,i}\subset U_{n,m}$
for all $i\le m$ we obtain
$$
\eta(X\backslash U_{n,m})\le
\nu_{1,m_n}(X\backslash U_{n,1})
+
\nu_{2,m_n}(X\backslash U_{n,2})
+\cdots
+
\nu_{m,m_n}(X\backslash U_{n,n-1})
\le 2^{-n}.
$$
Whenever $x\in U_{n,m}$ we have
$d(\xi_n^x,\xi_m^x)\le 2^{-n}$. Therefore,
\begin{align*}
&\biggl|
\int_X \psi(x)\int_Y \varphi(y)\, \xi_n^x(dy)\, \eta(dx)
-
\int_X \psi(x)
\int_Y \varphi(y)\, \xi_m^x(dy)\, \eta(dx)\biggr|
\\
&
\le \int_X\biggl|\int_Y \varphi(y)\, \xi_n^x(dy)
-\int_Y \varphi(y)\, \xi_m^x(dy)\biggr| |\psi(x)|\, \eta(dx)
\\
&
\le
\int_X
d(\xi_n^x,\xi_m^x)\, \eta(dx)
\le 2^{-n}+\eta(X\backslash U_{n,m})
\le 2^{1-n},
\end{align*}
whence (3) follows.

The continuity of the mapping $x\mapsto \xi_m^x$ with values
in $(\mathcal{P}(Y),d)$ yields the continuity of the function
$$
x\mapsto \int_Y \varphi(y)\, \xi_m^x (dy).
$$
Therefore, there exists $N>m$ such that
$$
\biggl|
\int_{X}\psi(x)\int_Y
\varphi(y)\, \xi_m^x (dy) \, \nu(dx)-
\int_{X}\psi(x)
\int_Y\varphi(y)\, \xi_m^x (dy) \, \nu_n(dx)
\biggr|
< \varepsilon \eqno(4)
$$
if  $n\ge N$. For such $n$ on account of relationships
(2), (3) and (4) we find
$$
\biggl|
\int_{X\times Y} f\, d\mu -\int_{X\times Y} f\, d\mu_n
\biggr|
\le 6\varepsilon .
$$
It remains to replace the measures $\mu_n$ by the probability measures
$\mu_n(X)^{-1}\mu_n$ and use the fact that
$\mu_n(X)=\nu(K_{n-1})\to 1$ as $n\to\infty$.
The theorem is proven.

\vskip .1in

{\bf Remark 2.}
(i)
It is obvious from the proof that our assumption on $Y$ can be weakened as follows:
it suffices that every Radon measure on $Y$ be concentrated on a countable union
of metrizable compact sets (i.e., be concentrated on a Souslin set).
For example, on can take for $Y$ a completely regular Souslin space.
It is easily verified that any continuous image of a compact with property (ud)
 has this property as well. More generally, properties (ud)
and (tud) are preserved by continuous mappings, for which the induced  mappings of the
 spaces of Radon probability  measures are surjective (for example,
 this is the case for continuous mappings, for which the preimages of
 compact sets are compact).

(ii)
Apparently, property (tud) is strictly stronger than property~(ud).
Cer\-tain\-ly, for compact spaces (more generally, for sequentially
Prohorov spaces) both properties are equivalent.
It is very likely that there exist compact spaces $X$ and $Y$ with property (ud),
whose product $X\times Y$ has no this property.
It would be interesting to construct the corresponding examples.

\vskip .1in

This work was supported by the RFBR projects 07-01-00536,
05-01-02941-JF, and GFEN-06-01-39003,
DFG Grant 436 RUS 113/343/0(R), ARC Discovery Grant DP0663153, and
the SFB 701 at the University of Bielefeld.
{\sloppy

}

\vskip .2in

\centerline{{\sc References}}

\vskip .1in

[1] Bogachev V.I. Measure theory. V.~2. Springer, Berlin -- New York,
2007 (2nd Russian ed.: Moscow, 2006).

[2]
Niederreiter H.
On the existence of uniformly distributed sequences
in compact spaces.
 Compositio Math. 1972. V.~25. P.~93--99.

[3]
Fremlin D.
Measure theory. V.~4.
University of Essex, Colchester, 2003.

[4]
Losert V.
On the existence of uniformly distributed sequences
in compact topological spaces.~ I.
 Trans. Amer. Math. Soc. 1978. V.~246. P.~463--471.

[5]
Losert V.
On the existence of uniformly
distributed sequences in compact topological spaces.~II.
 Monatsh. Math. 1979. B.~87, N~3. S.~247--260.

[6]
Mercourakis S.
Some remarks on countably determined
measures and uniform distribution of sequences.
Monatsh. Math. 1996. B.~121, N~1-2. S.~79--111.

[7]
Plebanek G.
Approximating Radon measures on first-countable compact spaces.
 Colloq. Math. 2000. V.~86, N~1. P.~15--23.

\vskip .2in

Department of Mechanics and Mathematics,

Moscow State University

119992 Moscow, Russia

}
\end{document}